\documentclass[11pt]{article}

\usepackage[T1]{fontenc}
\usepackage[utf8]{inputenc}
\usepackage{lmodern}

\usepackage{amsmath, amssymb, amsthm}
\usepackage{mathrsfs}

\usepackage[margin=2.5cm]{geometry}
\usepackage{microtype}

\usepackage{xcolor, graphicx, float}
\usepackage[breakable]{tcolorbox}
\usepackage[shortlabels]{enumitem}
\usepackage{tikz}
\usepackage{caption}
\usepackage{subcaption}
\usepackage{bbm}
\usepackage{soul}
\usepackage{empheq} 

\usepackage{etoolbox}

\makeatletter
\patchcmd{\@maketitle}{\LARGE \@title}{\LARGE\bfseries\@title}{}{}

\renewcommand{\@seccntformat}[1]{\csname the#1\endcsname.\quad}

\renewcommand\paragraph{\@startsection{paragraph}{4}{\z@}%
    {1.5ex \@plus1ex \@minus.2ex}%
    {-1em}%
    {\normalfont\normalsize\bfseries}}
\makeatother

\makeatletter
\def\th@plain{%
	\thm@notefont{}
	\itshape 
}
\def\th@definition{%
	\thm@notefont{}
	\normalfont 
}

\renewenvironment{proof}[1][\proofname]{\par
	\normalfont
	\topsep0\p@\@plus3\p@ \trivlist
	\item[\hskip\labelsep\itshape
	#1\@addpunct{.}]\ignorespaces
}{%
	\qed\endtrivlist
}
\makeatother

\newtheorem{theorem}{Theorem}[section]
\newtheorem{lemma}[theorem]{Lemma}
\newtheorem{corollary}[theorem]{Corollary}

\theoremstyle{definition}

\newtheorem{remark}[theorem]{Remark}
\newtheorem{assumption}[theorem]{Assumption}
\newtheorem{algorithm}{Algorithm}

\definecolor{darkblue}{rgb}{0,0,.5}

\usepackage{hyperref}
\hypersetup{
	colorlinks=true,		
	linkcolor=darkblue,		
	citecolor=darkblue,		
	urlcolor=darkblue 		
}

\allowdisplaybreaks

\newcommand{\dom}{\ensuremath{\operatorname{dom}}}
\newcommand{\Fix}{\ensuremath{\operatorname{Fix}}}
\newcommand{\zer}{\ensuremath{\operatorname{zer}}}
\newcommand{\gra}{\ensuremath{\operatorname{gra}}}
\newcommand{\Id}{\ensuremath{\operatorname{Id}}}

\newcommand{\ran}
{\ensuremath{\operatorname{ran}}}

\newcommand{\lspan}{\ensuremath{\operatorname{span}}}

\newcommand{\diag}{\ensuremath{\operatorname{diag}}}

\newcommand{\Lap}{\ensuremath{\operatorname{Lap}}}

\newcommand{\scal}[2]{\left\langle #1,#2 \right\rangle}

\def\bA{{\boldsymbol{A}}}
\def\bB{{\boldsymbol{B}}}
\def\bC{{\boldsymbol{C}}}

\def\bOne{{\mathbbm{1}}}

\def\bv{{\mathbf{v}}}

\def\bx{{\mathbf{x}}}
\def\by{{\mathbf{y}}}
\def\bz{{\mathbf{z}}}
\def\be{{\mathbf{e}}}

\begin{document}

\title{Weak and strong convergence of a general splitting approach for monotone multioperator inclusions}

\author{
Minh N. Dao\thanks{School of Science, RMIT University, Melbourne, VIC 3000, Australia.
E-mail:~\href{href:minh.dao@rmit.edu.au}{minh.dao@rmit.edu.au}.},
~
Matthew K. Tam\thanks{School of Mathematics and Statistics, The University of Melbourne, Parkville, VIC 3010, Australia.
E-mail:~\href{href:matthew.tam@unimelb.edu.au}{matthew.tam@unimelb.edu.au}.},
~and~
Thang D. Truong\thanks{School of Science, RMIT University, Melbourne, VIC 3000, Australia.
E-mail:~\href{href:thang.tdk64@gmail.com}{thang.tdk64@gmail.com}.}
}

\date{\today}

\maketitle

\begin{abstract}
In this paper, we study a general forward-backward splitting approach for monotone inclusion problems that simultaneously accommodate set-valued, monotone Lipschitz continuous, and cocoercive operators. Weak convergence is established  under a general setting. Strong convergence is obtained when both the maximally monotone operators are linear relations and the monotone Lipschitz continuous and cocoercive operators are linear. When the set-valued operators are normal cones to closed linear subspaces, we derive an explicit characterization of limit points of generated sequences using a new technique. Numerical experiments on a quadratic minimization problem with multiple-set split feasibility constraints are presented to illustrate the theoretical results.
\end{abstract}


\noindent{\bf Mathematics Subject Classification (MSC 2020):}
47J25,  
47J26,  
47H05,  
65K10.  

\section{Introduction}

We consider the inclusion problem
\begin{align}
\label{prob:ABC}
\text{find~} x \in \mathcal{H} \ \text{~such that~} \ 0 \in \sum_{i=1}^n A_i x + \sum_{k=1}^r B_k x + \sum_{j=1}^p C_j x,
\end{align}
where $\mathcal{H}$ is a real Hilbert space, $A_i\colon \mathcal{H}\rightrightarrows \mathcal{H}$ is maximally monotone for $i\in\{1,\dots,n\}$, $B_k\colon \mathcal{H}\to \mathcal{H}$ is monotone and Lipschitz continuous for $k\in\{1,\dots, r\}$, and $C_j\colon \mathcal{H}\to \mathcal{H}$ is cocoercive for $j\in\{1,\dots,p\}$. Problem \eqref{prob:ABC} encompasses optimization problems, saddle-point problems, and variational inequalities, with applications in signal processing and other areas \cite{BC17, BD18, CW05}.

When $r=0$ and $p=0$, available splitting methods include the \emph{proximal point algorithm} \cite{Roc76} for $n=1$, the \emph{Douglas--Rachford algorithm} \cite{DR56, Sva11} for $n=2$, and more recent algorithms \cite{BCN24, MT23, Ryu20, Tam23} for $n \geq 3$. These algorithms rely on evaluating the resolvents of the set-valued operators. Such evaluations, called \emph{backward steps}, may not admit closed-form expressions and can be computationally demanding. Therefore, when single-valued operators are present, \emph{forward steps} can exploit their direct evaluations instead. Specifically, when there are only cocoercive operators, i.e., $r=0$, the \emph{forward-backward algorithm} \cite{LM79, Passty79} applies to $n=p=1$, the \emph{Davis--Yin algorithm} \cite{DY17} to $n=2$ and $p=1$, while more general settings with $n \geq 2$ and $p\geq 1$ are considered in \cite{ACGN25, ACL25, AMTT23, DTT26}. When there is no cocoercive operators, i.e., $p=0$, the \emph{forward-backward-forward algorithm} \cite{Tseng00}, the \emph{forward-reflected-backward algorithm} \cite{MT20}, and the \emph{reflected-forward-backward algorithm} \cite{CV21} address the case $n =r =1$, whereas the algorithms in \cite{AMTT23, DTT26} accommodate $n\geq 3$ and $r\geq 1$. Although every cocoercive operator is monotone and Lipschitz continuous, the converse does not hold in general; hence, it is desirable to exploit the properties of different operators individually. For $n =r =p =1$, this motivates methods such as the \emph{forward-backward-half forward algorithm} \cite{BD18}, a variant of the forward-reflected-backward algorithm, and other approaches for monotone Lipschitz operators \cite{Rol25, ZCHD26}. 

The aforementioned algorithms generally guarantee only weak convergence of the generated sequences to a solution. Strong convergence typically requires additional assumptions on the involved operators such as strong monotonicity or alternative techniques \cite{BC01, BCM19, Xu02}. Linearity provides another setting in which strong convergence can be established. In particular, strong convergence of the Douglas--Rachford algorithm, the Ryu algorithm \cite{Ryu20}, and the Malitsky--Tam algorithm \cite{MT23} was established in \cite{BMSW25}. These results were subsequently extended in \cite{ABCL26} to a class of graph splitting methods for $n\geq 2$ and $r =p =0$. More recently, \cite{BBG26} established strong convergence under a divergent-series condition on the relaxation parameters. These developments motivate extending the framework to incorporate both monotone Lipschitz and cocoercive operators while retaining strong convergence. Explicit characterizations of the limit points can also provide additional information about the solutions selected by the algorithms. It is therefore desirable to derive such characterizations in special cases, for example, when the maximally monotone operators are normal cones of closed linear subspaces.

Our main contributions are as follows.
\begin{enumerate}
    \item 
    We study a general splitting approach for solving problem \eqref{prob:ABC} and establish weak convergence under general assumptions, with improved results when $B_1,\dots,B_r$ are linear, see Theorem~\ref{t:cvg_weak}. Strong convergence is guaranteed when the involved set-valued operators are linear relations and both monotone Lipschitz and cocoercive operators are linear, even when the relaxation parameter sequence only satisfies a divergent-series condition, see Theorem~\ref{t:cvg_strong}.
    \item 
    When the linear relations are specialized as normal cones to closed linear subspaces, we derive explicit expressions for the limit points in terms of projections onto the corresponding fixed point sets using a new technique, see Theorem~\ref{t:explicit}. The results are illustrated numerically on a quadratic minimization problem with multiple-set split feasibility constraints.
\end{enumerate}

The remainder of this paper is organized as follows. In Section~\ref{s:pre}, we introduce notation and background material on set-valued and single-valued operators. We present our main results about the convergence properties and limit points of generated sequences in Section~\ref{s:cvg}. Section~\ref{s:numerics} provides numerical experiments illustrating the computation of the limit points.

\section{Preliminaries}
\label{s:pre}

Throughout this paper, we write $\mathbb{N}$ and $\mathbb{R}$ for the set of nonnegative integers and real numbers. Let $\mathcal{H}$ be a real Hilbert space endowed with the inner product $\scal{\cdot}{\cdot}$ and induced norm $\|\cdot\|$. We use $\to$ and $\rightharpoonup$ to denote strong and weak convergence, respectively. 

The notation $A \colon \mathcal{H} \rightrightarrows \mathcal{H}$ indicates that $A$ is a set-valued operator, while $A \colon \mathcal{H}\to\mathcal{H}$ indicates that $A$ is a single-valued operator. Let $A$ be an operator on $\mathcal{H}$. The \emph{domain} and \emph{range} of $A$ are given by
\begin{align*}
    \dom A := \{x \in \mathcal{H} : Ax \neq \varnothing\}, 
    \quad \ran A := \{y \in \mathcal{H} : y \in Ax \text{~for some~} x \in \mathcal{H}\}.
\end{align*}
The \emph{graph} and \emph{set of fixed points} of $A$ are
\begin{align*}
    \gra A := \{(x, u) \in \mathcal{H} \times \mathcal{H} : u \in Ax\}, \quad  \Fix A := \{x \in \mathcal{H} : x \in Ax\}.
\end{align*}
The resolvent of $A\colon \mathcal{H} \to \mathcal{H}$ with parameter $\gamma > 0$ is $J_{\gamma A}:= (\Id + \gamma A)^{-1}$, where $\Id \colon \mathcal{H}\to\mathcal{H}$ denotes the identity operator. The operator $A$ is \emph{monotone} if, for all $(x, u), (y, v) \in \gra A$,
\begin{align*}
    \scal{x-y}{u-v} \geq 0
\end{align*}
and maximally monotone if its graph is not properly contained in the graph of any other monotone operator.

A single-valued operator $C \colon \mathcal{H}\to\mathcal{H}$ is \emph{$\ell$-Lipschitz continuous} for $\ell \in [0, +\infty)$ if, for all $x,y\in\mathcal{H}$,
\begin{align*}
    \|Cx - Cy \| \leq \ell \|x - y \|.
\end{align*}
When $\ell=1$, $C$ is \emph{nonexpansive}. An operator $C$ is \emph{$\frac{1}{\ell}$-cocoercive} for $\ell \in (0,+\infty)$ if, for all $x, y\in\mathcal{H}$,
\begin{align*}
    \scal{Cx - Cy}{x- y} \geq \frac{1}{\ell} \|Cx - Cy\|^2.
\end{align*}
Every $\frac{1}{\ell}$-cocoercive operator is monotone and $\ell$-Lipschitz continuous by the Cauchy--Schwarz inequality. In addition, every monotone and Lipschitz continuous operator is maximally monotone \cite[Corollary~20.28]{BC17}.

Let $T \colon \mathcal{H} \to \mathcal{H}$. We recall from \cite[Definition~2.1]{BDP22} that $T$ is said to be \emph{conically $\rho$-averaged} if $\rho \in (0, +\infty)$ and for all $x, y\in \mathcal{H}$, 
\begin{align*}
\|Tx -Ty\|^2 +\frac{1-\rho}{\rho}\|(\Id -T)x -(\Id -T)y\|^2 \leq \|x -y\|^2
\end{align*}
and from \cite[Definition~2.1]{DTT26} that $T$ is \emph{conically $\rho$-quasiaveraged} if $\rho\in (0, +\infty)$ and for all $x\in \mathcal{H}$ and all $y\in \Fix T$, 
\begin{align*}
\|Tx -y\|^2 +\frac{1-\rho}{\rho}\|(\Id -T)x\|^2 \leq \|x -y\|^2.
\end{align*}
It is clear from the definition that every conically $\rho$-averaged operator is conically $\rho$-quasiaveraged. Moreover, a conically $\rho$-averaged operator is nonexpansive when $\rho \in (0,1]$.

For a matrix $M$, we use $\ran M$ and $\ker M$ to denote its \emph{range} and \emph{kernel}, respectively, $M^\top$ for its transpose, and $M^\dagger$ for its Moore--Penrose pseudoinverse. We write $M \succeq 0$ to indicate that $M$ is positive semidefinite. We identify a matrix $M=(M_{ij})\in\mathbb{R}^{n \times m}$ with the Kronecker product $(M \otimes \Id) \colon \mathcal{H}^m \to \mathcal{H}^n$ whenever the context is clear. In particular, for any $\bz \in \mathcal{H}^m$, we write $(M\otimes \Id)\bz$ simply as $M\bz$. Similarly, for operators $\bB \colon \mathcal{H}^n \to \mathcal{H}^n$ and $\bC \colon \mathcal{H}^m \to \mathcal{H}^m$, we interpret the products $\bB M$ and $M\bC$ as $\bB(M\otimes \Id)$ and $(M\otimes \Id)\bC$, respectively. Finally, for any $\bx = (x_1,\dots,x_n) \in \mathcal{H}^n$, we define the diagonal operator $\bA \colon \mathcal{H}^n \rightrightarrows \mathcal{H}^n$ by $\bA \bx := A_1 x_1 \times \dots \times A_n x_n$.

\section{A general splitting approach and its convergence analysis}
\label{s:cvg}
Throughout this section, for each $i \in\{1,\dots,n\}$, $A_i$ is a maximally monotone operator, for each $k\in\{1,\dots,r\}$, $B_k$ is monotone and $\ell_k$-Lipschitz continuous with $\ell_k\in(0,+\infty)$, for each $j\in\{1,\dots p\}$, $C_j$ is $\frac{1}{\beta_j}$-cocoercive with $\beta_j \in (0,+\infty)$. Let $M =(M_{ij}) \in \mathbb{R}^{n \times m}$, $N =(N_{ij}) \in \mathbb{R}^{n \times n}$, $D=\diag(\delta_1,\dots,\delta_n)$ where $ \delta_1,\dots,\delta_n \in (0, +\infty)$, $P =(P_{ij})\in\mathbb{R}^{n\times r}$, $Q =(Q_{ij})\in\mathbb{R}^{n\times r}$, $R =(R_{ij}) \in \mathbb{R}^{r\times n}$, $H =(H_{ij})\in\mathbb{R}^{n\times p}$, and $K =(K_{ij}) \in \mathbb{R}^{p\times n}$. Define an auxiliary operator $\Phi :\mathcal{H}^n\to\mathcal{H}^n$ by
\begin{align}\label{e:Phi}
\Phi := (P-Q)\bB R + Q\bB P^\top + H\bC K
\end{align}
We denote $\bOne_n=(1,\dots,1)\in\mathbb{R}^n$ and omit the subscript when the context is clear. Let $Z :=\zer\left(\sum_{i=1}^n A_i +\sum_{k=1}^r B_k +\sum_{j=1}^p C_j\right)$ be the solution set of problem \eqref{prob:ABC}. We denote $\Delta := \{\bx = (x_1,\dots,x_n)\in \mathcal{H}^n \ |\ x_1 = \dots = x_n\}$, $\diag(\ell) :=\diag(\ell_1,\dots,\ell_r)$, $\diag(\beta) := \diag(\beta_1,\dots,\beta_p)$, and set
\begin{align*}
\Omega &:= 2D - N - N^\top - MM^\top, \\
\Upsilon &:= (P-R^\top)\diag(\ell)(P^\top -R) + (P-Q)\diag(\ell)(P^\top -Q^\top)  + \frac{1}{2}(H-K^\top)\diag(\beta)(H^\top - K),
\end{align*}
and consider the following assumptions on the coefficient matrices.

\begin{assumption}[Standing assumptions]
\label{a:stand}
~
\begin{enumerate}
 \item\label{a:stand_kerM}
$\ker M^\top = \lspan\{\bOne\}$.
\item\label{a:stand_N}
$\sum_{i,j=1}^n N_{ij} = \bOne^\top N \bOne= \bOne^\top D \bOne= \sum_{i=1}^n \delta_i$.
\item\label{a:stand_PQR}
$P^\top \bOne= Q^\top \bOne=R\bOne=\bOne$.
\item \label{a:stand_HK}
$H^\top \bOne = K \bOne = \bOne$.
\end{enumerate}
\end{assumption}

\begin{assumption}[Explicit assumptions]
\label{a:explicit}
~
\begin{enumerate}
\item 
$N_{ij}=0$ whenever $1 \leq i \leq j \leq n$.
\item 
$(P_{ik} - Q_{ik}) R_{kl} = 0$ and $Q_{ik}P_{lk}=0$ whenever $l\geq i$ and $k \in \{1,\dots,r\}$.
\item 
$H_{ij}K_{jl}=0$ whenever $l\geq i$ and $j\in\{1,\dots,p\}$.
\end{enumerate}
\end{assumption}

\begin{assumption}[Positive semidefiniteness]  
\label{a:PSD} 
There exists $\alpha \in [0,1)$ and $\gamma \in (0,+\infty)$ such that
\begin{align*}
\Omega + \alpha MM^\top -\gamma \Upsilon \succeq 0.
\end{align*}
\end{assumption}

\begin{remark}[Implications of Assumptions~\ref{a:stand} and \ref{a:PSD}]
\label{r:compute_stepsize}
Under Assumption~\ref{a:stand}, by the same arguments as in \cite[Remark~3.3(iii)]{DTT26}, $U = (P^\top - R)(M^\top)^\dagger \in \mathbb{R}^{r \times m}$ solves $UM^\top = P^\top - R$ with minimal norm, $V = (P^\top - Q^\top)(M^\top)^\dagger \in \mathbb{R}^{r \times m}$ solves $VM^\top = P^\top - Q^\top$ with minimal norm, and $X = (H^\top - K)(M^\top)^\dagger \in \mathbb{R}^{p \times m}$ solves $XM^\top = H^\top - K$ with minimal norm. Let 
\begin{align}
\label{def:Theta}
\Theta := U^\top \diag(\ell)U + V^\top \diag(\ell) V + \frac{1}{2}X^\top \diag(\beta) X.
\end{align}
Assumption~\ref{a:PSD} becomes 
\begin{align}
\label{r:PSD}
\Omega + M(\alpha \Id - \gamma \Theta)M^\top \succeq 0. 
\end{align}
\end{remark}

\begin{remark}[Comparing with existing works]
\label{r:compare_stepsize}
For $\alpha=0$, Assumption~\ref{a:PSD} becomes
\begin{align*}
\Omega - \gamma \Upsilon \succeq 0,    
\end{align*}
which coincides with \cite[Assumption~3.1(c)]{ZCHD26}. When $\bC=0$, this condition reduces to \cite[Remark~3.7(i)]{DTT26} and when $\bB = 0$, it becomes \cite[Assumption~4.7]{ACGN25} and \cite[Remark~3.7(ii)]{DTT26}.
    
For $\alpha \in (0,1)$, one sufficient condition for \eqref{r:PSD} is 
\begin{align*}
\Omega \succeq 0 \text{~~and~~} \gamma \leq \frac{\alpha}{\lambda_{\max}(\Theta)}.
\end{align*}
The latter is satisfied as soon as $\gamma < \frac{1}{\lambda_{\max}(\Theta)}$ and $\alpha = \gamma \lambda_{\max}(\Theta)$, see \cite[Remark~3.11(i)]{DTT26}. There is a broad class of splitting algorithms that allow choosing $M$, $N$, and $D$ based on graphs and satisfy the condition $\Omega \succeq 0$ automatically. A special case is when $\Omega = \kappa MM^\top$, condition \eqref{r:PSD} reduces to 
\begin{align*}
M(\kappa + \alpha - \gamma \Theta )M^\top \succeq 0,    
\end{align*}
yielding a larger upper bound $\gamma \leq \frac{\kappa + \alpha}{\lambda_{\max}(\Theta)}$, see \cite[Section~4]{DTT26} and \cite[Example~4.2]{DPTT26}.  
\end{remark}

An approach for solving problem \eqref{prob:ABC} when either $\bB=0$ or $\bC=0$, together with its convergence analysis was proposed in \cite[Theorem~3.9 and 3.10]{DTT26}. The following algorithm is a natural extension of \cite[Algorithm~1]{DTT26}, and its weak convergence under the restriction $\alpha=0$ in Assumption~\ref{a:PSD} was subsequently established in \cite[Theorem~4.1]{ZCHD26}. We provide a proof of weak convergence under general assumptions in Theorem~\ref{t:cvg_weak} and strong convergence when the involved set-valued operators are linear relations and both monotone Lipschitz and cocoercive operators are linear in Theorem~\ref{t:cvg_strong}.
\begin{tcolorbox}
[boxsep=5pt,left=0pt,right=0pt,top=0pt,bottom=0pt, colback=white,colframe=black,
    boxrule=0.5pt,breakable]
\begin{algorithm}
\label{algo:full}
Let $\gamma \in (0, +\infty)$ and $(\lambda_t)_{t\in\mathbb{N}}\subset (0, +\infty)$. Let $\bz^0 = (z^0_1, \dots, z^0_{m})\in\mathcal{H}^{m}$. For each $t\in\mathbb{N}$, compute
\begin{align*}
\begin{cases}
x^t_1 &= J_{\frac{\gamma}{\delta_1}A_1}\Big(\frac{1}{\delta_1}\sum_{j=1}^{m} M_{1j} z_j^t \Big)\\
x^t_i &=J_{\frac{\gamma}{\delta_i}A_i}\Bigg(\frac{1}{\delta_i}\Big(\sum_{j=1}^{m} M_{ij}z^t_j +\sum_{j=1}^{i-1} N_{ij}x^t_j - \gamma \sum_{j=1}^r (P_{ij} -Q_{ij})B_j\left(\sum_{l=1}^{i-1} R_{jl}x^t_l\right) \\ 
&\qquad\qquad - \gamma \sum_{j=1}^r Q_{ij}B_j\left(\sum_{l=1}^{i-1} P_{lj}x^t_l\right)\Big) - \gamma \sum_{j=1}^p H_{ij}C_j\left(\sum_{l=1}^{i-1} K_{jl} x_l^t\right) \Bigg),\quad i \in \{2, \dots, n\} \\
z^{t+1}_i \!\!\!\!\! &=z^t_i -\lambda_t \sum_{j=1}^n M_{ji}x^t_j,\quad i\in \{1, \dots, m\},
\end{cases}
\end{align*}
or compactly as
\begin{align*}
\begin{cases}
\bx^t &=J_{\gamma D^{-1}\bA}( D^{-1}( M\bz^t +  N\bx^t -\gamma\Phi \bx^t)) \\
\bz^{t+1}\!\!\!\!\! &=\bz^t - \lambda_t M^{\top}\bx^t.
\end{cases}
\end{align*}
\end{algorithm}
\end{tcolorbox}

We note that the sequence $(\bz^t)_{t\in\mathbb{N}}$ generated by Algorithm~\ref{algo:full} can be obtained via
\begin{align*}
\bz^{t+1} = \left(1 - \lambda_t\right)\bz^t + \lambda_tT\bz^t,
\end{align*}
where $T\colon \mathcal{H}^{m} \to \mathcal{H}^{m}$ is defined by 
\begin{align}
\label{def:T}
T\bz = \bz - M^\top \bx
\end{align}
with $\bx \in \mathcal{H}^n$ given via the solution operator $S\colon \mathcal{H}^{m} \to \mathcal{H}^n$
\begin{align}
\label{def:S}
\bx = S\bz := J_{\gamma D^{-1}\bA}( D^{-1}( M\bz +  N\bx -\gamma\Phi\bx)).
\end{align}
We introduce the vector 
\begin{align}
\label{def:mu}
\boldsymbol{\eta} :=M^\dagger (D -N)\bOne \in \mathbb{R}^m.    
\end{align}
By Assumption~\ref{a:stand}\ref{a:stand_kerM}--\ref{a:stand_N}, $(D -N)\bOne \in \ran M$, and so $M\boldsymbol{\eta} =(D-N)\bOne$. Note that, since the first row of $N$ is zero and the first row of $D$ is $\delta_1 > 0$, $(D-N)\bOne \neq 0$, and hence $\boldsymbol{\eta} \neq 0$.

\begin{lemma}[Fixed points] 
\label{l:fixed_points}
Suppose that Assumption~\ref{a:stand} holds. Then $S(\Fix T) =\{\bOne x: x\in Z\}$. Consequently, $\Fix T \neq \varnothing$ if and only if $Z \neq \varnothing$. Moreover, for all $\bz\in \mathcal{H}^m$ and $x\in \mathcal{H}$, 
\begin{align*}
S\bz = \bOne x \iff M(\bz -\boldsymbol{\eta}x)\in \gamma(\bA +P\bB +H\bC)(\bOne x).
\end{align*}
\end{lemma}
\begin{proof}
By \cite[Lemma~4.1]{ZCHD26} (see also \cite[Lemma~3.5]{DTT26}), every $\bz\in \Fix T$ satisfies $S\bz =\bOne x$ for some $x\in Z$, and every $x\in Z$ admits such a fixed point $\bz$. The first two assertions follow. 

Next, let $\bz\in\mathcal{H}^m$ and $x\in\mathcal{H}$. By the definition of $S$, for all $\bx\in \mathcal{H}^n$,
\begin{align*}
S\bz =\bx 
\iff M\bz +N\bx -\gamma\Phi\bx\in D\bx +\gamma\bA\bx 
\iff M\bz -(D -N)\bx\in \gamma(\bA +\Phi)\bx.
\end{align*}
It follows from Assumption~\ref{a:stand}\ref{a:stand_PQR}--\ref{a:stand_HK} and \eqref{def:mu} that $\Phi(\bOne x) =P\bB(\bOne x) +H\bC(\bOne x)$ and $(D -N)\bOne x =M\boldsymbol{\eta}x$. Substituting $\bx =\bOne x$ into the preceding equivalence yields the final assertion.
\end{proof}

\begin{lemma}[Quasiaveragedness and averagedness]
\label{l:quasiaveraged}
Suppose that Assumption~\ref{a:stand} holds. Then the following hold:
\begin{enumerate}
\item \label{averaged:ineq}
For all $\bz \in \mathcal{H}^m$ and all $\bz' \in \Fix T$, 
\begin{multline}
\label{ineq:T}
\|T\bz - T\bz'\|^2 -\alpha \|(\Id - T)\bz - (\Id -T)\bz'\|^2 \\ +\scal{\bx - \bx'}{(\Omega + \alpha MM^\top - \gamma\Upsilon)(\bx - \bx')} \leq \|\bz - \bz'\|^2,
\end{multline}
where $\bx =S\bz$ and $\bx' =S\bz'$. Moreover, if $B_1,\dots,B_r$ are linear, then this inequality holds for all $\bz'\in \mathcal{H}^m$.
\item \label{aveaged:T}
If Assumption~\ref{a:PSD} holds, then $T$ is conically $\frac{1}{1-\alpha}$-quasiaveraged and, if additionally $B_1,\dots,B_r$ are linear, then $T$ is conically $\frac{1}{1-\alpha}$-averaged.
\end{enumerate}
\end{lemma}
\begin{proof}
\ref{averaged:ineq} \& \ref{aveaged:T}: 
Let $\bz, \bz' \in \mathcal{H}^m$, $\bx = S\bz$, and $\bx' = S\bz'$. Applying \cite[Lemma~3.8]{DPTT26} with $\Gamma = \Id -T$ and noting that its proof does not use any assumptions on $\Phi$, we obtain that
\begin{multline*}
    \scal{(\Id - T)\bz - (\Id-T)\bz'}{\bz - \bz'} \geq \frac{1-\alpha}{2}\|(\Id - T)\bz - (\Id-T)\bz'\|^2 + \gamma \scal{\bx - \bx'}{\Phi\bx - \Phi \bx'} \\
    + \frac{1}{2}\scal{\bx - \bx'}{(\Omega + \alpha MM^\top)(\bx - \bx')},
\end{multline*}
which is equivalent to
\begin{multline}
\label{eq:T_Phi}
    \|T\bz - T\bz'\|^2  -\alpha \|(\Id - T)\bz - (\Id -T)\bz'\|^2 + \scal{\bx - \bx'}{(\Omega + \alpha MM^\top)(\bx - \bx')} \\
    + 2\gamma \scal{\bx - \bx'}{\Phi\bx - \Phi\bx'} \leq \|\bz - \bz'\|^2.
\end{multline}
Recalling the definition of $\Phi$ in \eqref{e:Phi}, we have 
\begin{align*}
\scal{\bx - \bx'}{\Phi\bx - \Phi\bx'} = \scal{\bx - \bx'}{\Phi_{\bB}\bx - \Phi_{\bB}\bx'} +  \scal{\bx - \bx'}{\Phi_{\bC}\bx - \Phi_{\bC}\bx'}    
\end{align*}
where $\Phi_{\bB} :=(P-Q)\bB R + Q\bB P^\top$ and $\Phi_{\bC} :=H\bC K$. It follows from the proof of \cite[Lemma~3.10]{DPTT26} that, for all $\bx, \bx' \in \mathcal{H}^n$,
\begin{align*}
\scal{\bx - \bx'}{\Phi_{\bC}\bx - \Phi_{\bC}\bx'} \geq -\frac{1}{4} \scal{\bx - \bx'}{\big((H-K^\top)\diag(\beta)(H^\top - K)\big)(\bx - \bx')}.
\end{align*}
Next, the proof of \cite[Lemma~3.9]{DPTT26} implies that, for all $\bx, \bx' \in \mathcal{H}^n$ with $P^\top\bx' = R\bx' = \bx'$, 
\begin{multline}
\label{eq:Phi_B}
\scal{\bx - \bx'}{\Phi_{\bB}\bx - \Phi_{\bB}\bx'}\\ \geq -\frac{1}{2} \scal{\bx - \bx'}{\big((P-R^\top)\diag(\ell)(P^\top - R) + (P-Q)\diag(\ell)(P^\top - Q^\top)\big)(\bx - \bx')}.
\end{multline}
By combining these inequalities, for all $\bx, \bx' \in \mathcal{H}^n$ with $P^\top\bx' = R\bx' = \bx'$, 
\begin{align}
\label{eq:Phi}
\scal{\bx - \bx'}{\Phi\bx - \Phi\bx'} \geq -\frac{1}{2}\scal{\bx -\bx'}{\Upsilon(\bx - \bx')}.
\end{align}
Note that $\bz' \in \Fix T$ implies that $P^\top\bx' = R\bx' = \bx'$. Therefore, substituting \eqref{eq:Phi} into \eqref{eq:T_Phi}, we obtain \eqref{ineq:T}, which means $T$ is conically $\frac{1}{1-\alpha}$-quasiaveraged when Assumption~\ref{a:PSD} holds.

Finally, applying \eqref{eq:Phi_B} with $\bx' =0$, we have for all $x\in \mathcal{H}^n$ that
\begin{align}
\label{eq:Phi_B'}
\scal{\bx}{\Phi_{\bB}\bx} \geq -\frac{1}{2} \scal{\bx}{\big((P-R^\top)\diag(\ell)(P^\top - R) + (P-Q)\diag(\ell)(P^\top - Q^\top)\big)\bx}.
\end{align}
If $B_1,\dots,B_r$ are linear, then replacing $\bx$ in \eqref{eq:Phi_B'} by $\bx -\bx'$, we obtain that \eqref{eq:Phi_B} holds for all $\bx, \bx' \in \mathcal{H}^n$, and so does \eqref{eq:Phi}. Therefore, \eqref{ineq:T} holds for all $\bz, \bz' \in \mathcal{H}^m$, which implies that $T$ is conically $\frac{1}{1-\alpha}$-averaged under Assumption~\ref{a:PSD}.
\end{proof}

\begin{theorem}[Weak convergence]
\label{t:cvg_weak}
Suppose that $Z \neq \varnothing$ and that Assumptions~\ref{a:stand}, \ref{a:explicit}, and \ref{a:PSD} hold. Let $(\bx^t)_{t\in\mathbb{N}}$ and $(\bz^t)_{t\in\mathbb{N}}$ be the sequences generated by Algorithm~\ref{algo:full} with an initial point $\bz^0\in\mathcal{H}^{m}$ and a sequence of relaxation parameters $(\lambda_t)_{t\in\mathbb{N}}$ in $[0,1-\alpha]$. Suppose further that either $\lim\inf_{t\to+\infty} \lambda_t(1-\frac{\lambda_t}{1-\alpha}) > 0$ or $B_1,\dots,B_r$ are linear and $\sum_{t=1}^{+\infty} \lambda_t(1-\frac{\lambda_t}{1-\alpha}) = +\infty$. Then the following hold:
\begin{enumerate}
    \item \label{t:cvg_weak_rate}
    $(\Id-T)\bz^t \to 0$ and $\sum_{i=1}^n s_i x_i \to 0$ whenever $(s_1,\dots,s_n) \in \mathbb{R}^n$ with $\sum_{i=1}^n s_i = 0$. Moreover, if $\lim\inf_{t\to+\infty} \lambda_t(1-\frac{\lambda_t}{1-\alpha}) > 0$, then $\|\frac{1}{t+1}\sum_{k=0}^t (\Id -T)\bz^k\| = O\big(\frac{1}{\sqrt{t}}\big)$ and, if additionally $B_1,\dots,B_r$ are linear, then $\|(\Id - T)\bz^t\| = o\big(\frac{1}{\sqrt{t}}\big)$.
    \item \label{t:cvg_weak_zx}
    $\bz^t \rightharpoonup \bar{\bz} \in \Fix T$ and $\bx^t \rightharpoonup \bar{\bx} = \bOne \bar{x}$ with $\bar{x} = J_{\frac{\gamma}{\delta_1}A_1}\Big(\frac{1}{\delta_1}\sum_{j=1}^{m} M_{1j} \bar{z}_j \Big) \in Z$.
\end{enumerate}
\end{theorem}
\begin{proof}
Since $Z \neq \varnothing$, it follows from Lemma~\ref{l:fixed_points} that $\Fix T \neq \varnothing$. In view of Lemma~\ref{l:quasiaveraged}, $T$ is conically $\frac{1}{1-\alpha}$-quasiaveraged and, if additionally $B_1,\dots,B_r$ are linear, then $T$ is conically $\frac{1}{1-\alpha}$-averaged. The conclusions follow by combining \cite[Proposition~2.2]{DTT26} and \cite[Proposition~2.9]{BDP22} with \cite[Lemma~3.8]{DTT26} and its proof.
\end{proof}

\begin{theorem}[Strong convergence]
\label{t:cvg_strong}
    Suppose that for each $i\in\{1,\dots,n\}$, $A_i$ is a linear relation, for each $k\in\{1,\dots,r\}$ and each $j\in\{1,\dots,p\}$, $B_k$ and $C_j$ are linear operators, that $Z \neq \varnothing$, and that Assumptions~\ref{a:stand}, \ref{a:explicit}, and \ref{a:PSD} hold. Let $(\bx^t)_{t\in\mathbb{N}}$ and $(\bz^t)_{t\in\mathbb{N}}$ be the sequences generated by Algorithm~\ref{algo:full} with an initial point $\bz^0\in\mathcal{H}^{m}$ and a sequence of relaxation parameters $(\lambda_t)_{t\in\mathbb{N}}$ in $[0,1-\alpha]$ satisfying $\sum_{t=1}^{+\infty} \lambda_t(1-\frac{\lambda_t}{1-\alpha}) = +\infty$. Then the following hold:
    \begin{enumerate}
        \item \label{t:ST}
        $S$ and $T$ are linear bounded operators. 
        \item \label{t:cvg_strong_zt}
        $\bz^t \to \bar{\bz} = P_{\Fix T}(\bz^0)$ and $\bx^t \to \bar{\bx} = \bOne \bar{x}$ with $\bar{x} = J_{\frac{\gamma}{\delta_1}A_1}\Big(\frac{1}{\delta_1}\sum_{j=1}^{m} M_{1j} \bar{z}_j \Big) \in Z$.
    \end{enumerate}  
\end{theorem}
\begin{proof}
\ref{t:ST}: Since $A_i$ is a linear relation and $\frac{\gamma}{\delta_i} > 0$, \cite[Theorem~2.1(xviii)]{BMW12} implies that $J_{\frac{\gamma}{\delta_i} A_i}$ is a linear operator. As $J_{\frac{\gamma}{\delta_i} A_i}$ is firmly nonexpansive, it is bounded. Moreover, since $B_k$ and $C_j$ are also bounded linear, using the definition of $S$, Assumption~\ref{a:stand} and \ref{a:explicit}, we can apply the linearity and boundedness inductively for $x_i^t$ and deduce that $S$ is bounded linear. Consequently, by the definition \eqref{def:T} and \eqref{def:S}, $T = \Id - M^\top S$ is also a bounded linear operator.  

\ref{t:cvg_strong_zt}: Since $B_1,\dots,B_r$ are linear, by Lemma~\ref{l:quasiaveraged}, $T$ is conically $\frac{1}{1-\alpha}$-averaged. We note that 
\begin{align*}
\bz^{t+1} =(1 -\lambda_t)\bz^t +\lambda_t T\bz^t =\left(1 -\frac{\lambda_t}{1-\alpha}\right)\bz^t +\frac{\lambda_t}{1-\alpha} T_{1-\alpha}\bz^t,    
\end{align*}
where $T_{1-\alpha} := \alpha\Id + (1-\alpha)T$. Since $T$ is linear, so is $T_{1-\alpha}$. Moreover, by \cite[Proposition~2.2(i) \& (ii)]{BDP22}, $T_{1-\alpha}$ is conically $1$-averaged, i.e., nonexpansive. Now, by applying \cite[Theorem~2.3]{BBG26}, $\bz^t \to \bar{\bz} = P_{\Fix T}(\bz^0)$. Since $S$ is continuous and $\bx^t = S\bz^t$, we deduce that $\bx^t \to \bar{\bx}$. In view of Theorem~\ref{t:cvg_weak}\ref{t:cvg_weak_zx}, $\bar{\bx} = \bOne \bar{x}$ with $\bar{x} = J_{\frac{\gamma}{\delta_1}A_1}\Big(\frac{1}{\delta_1}\sum_{j=1}^{m} M_{1j} \bar{z}_j \Big) \in Z$.
\end{proof}

\begin{remark}[Some instances of Algorithm~\ref{algo:full}]
\label{r:algorithms}
We provide some concrete instances of Algorithm~\ref{algo:full} for solving problem \eqref{prob:ABC} when $r = n-2$ and $p=m=n-1$. Let $G$ be a graph with uniform edge weight $w > 0$ and let $G'$ be an oriented subgraph $G'$ with uniform edge weights $\mu^2 \in (0, w]$. The matrices $N$ and $D$ are defined as the strictly lower triangular part of the \emph{weight matrix} and half of the \emph{degree matrix} of $G$, respectively, while $M$ is defined such that $MM^\top = \Lap(G')$ with $\Lap(G')$ as the \emph{Laplacian matrix} of $G'$. For example, the complete-parallel down algorithm means $G$ is a complete graph, whereas $G'$ is a parallel down graph (star graph with node $n$ as the central node). We refer readers to \cite[Section~4]{DTT26} for more details. For simplicity, we use the following $P$, $Q$, $R$, $H$, and $K$ for all considered algorithms 
\begin{align*}
      P &=\left[ 
        \begin{array}{c} 
          0_{1\times (n-2)} \\ 
          \hline 
          \Id_{n-2} \\
          \hline
          0_{1\times (n-2)}
        \end{array} 
        \right],
    &
    Q &=\left[ 
        \begin{array}{c} 
          0_{2\times (n-2)} \\ 
          \hline 
          \Id_{n-2} 
        \end{array} 
        \right], 
  & 
    R &=\left[ 
        \begin{array}{c|c} 
          \Id_{n-2} & 0_{(n-2)\times 2}
        \end{array} 
        \right], \\
    H &=\left[ 
        \begin{array}{c} 
          0_{1\times (n-1)} \\ 
          \hline 
          \Id_{n-1}
        \end{array} 
        \right], 
  & &
  & 
    K &=\left[ 
        \begin{array}{c|c} 
          \Id_{n-1} & 0_{(n-1)\times 1}
        \end{array} 
        \right],
\end{align*}
which satisfy Assumption~\ref{a:stand}. Then $\Phi\bx^t = \by^t = (y_1^t,y_2^t,\dots,y_n^t)$, where
\begin{align*}
    \begin{cases}
        y_1^t &= 0\\
        y_2^t &= B_1 x_1^t + C_1 x_1^t\\
        y_i^t &= (B_{i-1} x_{i-1}^t - B_{i-2}x_{i-2}^t) + B_{i-2} x_{i-1}^t + C_{i-1} x_{i-1}^t, \quad i\in\{3,\dots,n-1\}\\
        y_n^t &= -B_{n-2}x_{n-2}^t + B_{n-2} x_{n-1}^t + C_{n-1} x_{n-1}^t.
    \end{cases}
\end{align*}

We obtain the following coefficient matrices and algorithms based on their corresponding graphs.
\begin{enumerate}[label=\arabic*.]
\item
Complete-parallel down algorithm
\begin{align}
    N_{ij} &= 
    \begin{cases}
        w &\text{~if~} i > j\\
        0 &\text{~otherwise}
    \end{cases},
    \quad
    D = \frac{w(n-1)}{2}\diag(1,\dots,1), \notag\\
    M_{ij} &= 
    \begin{cases}
        \mu  &\text{~if~} i = j\\
        -\mu  &\text{~if~} i=n-1, j \in\{1,\dots, n-1\}\\
        0 &\text{~otherwise},
    \end{cases} \notag\\
&\begin{cases}
x_1^t &= J_{\frac{2\gamma}{w(n-1)}A_1}\left(\frac{2}{w(n-1)} \mu z_1^t\right)\\
x_i^t &= J_{\frac{2\gamma}{w(n-1)}A_i}\left(\frac{2}{w(n-1)}\left(\sum_{j=1}^{i-1} w x_j^t + \mu z_i^t - \gamma y_i^t\right)\right),\quad i\in\{2,\dots, n-1\}\\
x_n^t &= J_{\frac{2\gamma}{w(n-1)}A_n}\left(\frac{2}{w(n-1)}\left(\sum_{j=1}^{n-1} w x_j^t - \sum_{j=1}^{n-1}\mu z_j^t - \gamma y_n^t \right)\right)\\
z_i^{t+1} &= z_i^t - \lambda_t \mu (x_i^t - x_n^t), \quad i\in\{1,\dots,n-1\},
\end{cases}
\end{align} 
which simplifies to \cite[Equation~(56)]{DTT26} when either $\bB =0$ or $\bC=0$.

\item
Ring-sequential algorithm
\begin{align}
    N &= 
    \begin{cases}
        w &\text{~if~} i = j+1 \text{~or~} i=n, j=1\\
        0 &\text{~otherwise}
    \end{cases},
    \
    D = \frac{w}{2}\diag(1,\dots,1),\
    M = 
    \begin{cases}
        \mu &\text{~if~} i = j\\
        -\mu &\text{~if~} i=j+1\\
        0 &\text{~otherwise},
    \end{cases}\notag\\
&\begin{cases}
x_1^t &= J_{\frac{\gamma}{w} A_1}(\frac{\mu}{w} z_1^t)\\
x_i^t &= J_{\frac{\gamma}{w} A_i}(\frac{1}{w}(wx_{i-1}^t + \mu z_i^t - \mu z_{i-1}^t - \gamma y_i)),\quad i\in\{2,\dots, n-1\}\\
x_n^t &= J_{\frac{\gamma}{w} A_n}(\frac{1}{w}(w x_1^t + w x_{n-1}^t - \mu z_{n-1}^t -\gamma y_n))\\
z_i^{t+1} &= z_i^t - \lambda_t \mu(x_i^t - x_{i+1}^t), \quad i\in\{1,\dots,n-1\},
\end{cases}
\end{align}
which reduces to \cite[Equations~(43) and (44)]{DTT26} when $(\bB=0, \bC\neq0)$ and $(\bB\neq 0, \bC=0)$, respectively.

\item
Parallel up algorithm
\begin{align}
    N &= 
    \begin{cases}
        w &\text{~if~} i\in\{2,\dots,n\}, j=1\\
        0 &\text{~otherwise}
    \end{cases},
    \quad
    D = \frac{w}{2}\diag \left(n-1,1,\dots,1\right), \notag\\
    M &= 
    \begin{cases}
        \mu &\text{~if~} i = 1, j \in\{1,\dots,n-1\}\\
        -\mu &\text{~if~} i=j+1\\
        0 &\text{~otherwise},
    \end{cases}\notag\\
&\begin{cases}
x_1^t &= J_{\frac{2\gamma}{w(n-1)}A_1}\left(\frac{2\mu}{w(n-1)}\sum_{j=1}^{n-1} z_j^t\right)\\
x_i^t &= J_{\frac{2\gamma}{w} A_i}(\frac{2}{w}(w x_1^t - \mu z_{i-1}^t - \gamma y_i)),\quad i\in\{2,\dots, n\}\\
z_i^{t+1} &= z_i^t - \lambda_t \mu(x_1^t - x_{i+1}^t), \quad i\in\{1,\dots,n-1\},
\end{cases}
\end{align}
which is similar to \cite[Equations~(45) and (46)]{DTT26} in the cases $(\bB=0, \bC\neq0)$ and $(\bB\neq 0, \bC=0)$, respectively, but differs in the choices of $H$, $K$ for the former and $P$, $Q$, $R$ for the latter.

\item
Parallel down algorithm
\begin{align}
    N &= 
    \begin{cases}
        w &\text{~if~} i=n, j\in\{1,\dots,n-1\}\\
        0 &\text{~otherwise}
    \end{cases},
    \quad
    D = \frac{w}{2}\diag\left(1,\dots,1, n-1\right), \notag\\
    M &= 
    \begin{cases}
        \mu &\text{~if~} i = j\\
        -\mu &\text{~if~} i=n, j\in\{1,\dots,n-1\}\\
        0 &\text{~otherwise},
    \end{cases}\notag \\
&\begin{cases}
x_1^t &= J_{\frac{2\gamma}{w} A_1}(\frac{2\mu}{w}z_1^t)\\
x_i^t &= J_{\frac{2\gamma}{w} A_i}(\frac{2}{w}(\mu z_i^t - \gamma y_i)),\quad i\in\{2,\dots, n-1\}\\
x_n^t &= J_{\frac{2\gamma}{w(n-1)} A_n}\left(\frac{2}{w(n-1)}\left(w\sum_{j=1}^{n-1}x_j^t - \mu\sum_{j=1}^{n-1}z_j^t - \gamma y_n\right)\right)\\
z_i^{t+1} &= z_i^t - \lambda_t \mu(x_i^t - x_n^t), \quad i\in\{1,\dots,n-1\}.
\end{cases}
\end{align}
which is closely related to \cite[Equations~(48) and (49)]{DTT26} in the cases $(\bB=0, \bC\neq0)$ and $(\bB\neq 0, \bC=0)$, respectively, but uses alternative choices of $H$, $K$ in the former case and $P$, $Q$, $R$ in the latter.

\item
Sequential algorithm
\begin{align}
\label{algo:seq}
    N &= 
    \begin{cases}
        w &\text{~if~} i = j+1\\
        0 &\text{~otherwise}
    \end{cases},
    \quad
    D = w\diag\left(\frac{1}{2},1,\dots,1, \frac{1}{2}\right),\quad
    M = 
    \begin{cases}
        \mu &\text{~if~} i = j\\
        -\mu &\text{~if~} i=j+1\\
        0 &\text{~otherwise},
    \end{cases}\notag\\
&\begin{cases}
x_1^t &= J_{\frac{2\gamma}{w} A_1}(\frac{2\mu}{w}z_1^t)\\
x_i^t &= J_{\frac{\gamma}{w} A_i}(\frac{1}{w}(wx_{i-1}^t + \mu z_i^t - \mu z_{i-1}^t - \gamma y_i)),\quad i\in\{2,\dots, n-1\}\\
x_n^t &= J_{\frac{2\gamma}{w} A_n}(\frac{2}{w}(w x_{n-1}^t - \mu z_{n-1}^t - \gamma y_n))\\
z_i^{t+1} &= z_i^t - \lambda_t \mu(x_i^t - x_{i+1}^t), \quad i\in\{1,\dots,n-1\},
\end{cases}
\end{align}
which becomes \cite[Equation~(49)]{DPTT26} without dual terms when $\bB=0$.

\item
Complete algorithm 
\begin{align}
    N &= 
    \begin{cases}
        w &\text{~if~} i > j\\
        0 &\text{~otherwise}
    \end{cases},
    \quad
    D = \frac{w(n-1)}{2}\diag(1,\dots,1),
    \quad
    M = 
    \begin{cases}
        a_i &\text{~if~} i = j\\
        -b_j &\text{~if~} i > j\\
        0 &\text{~otherwise},
    \end{cases}\notag\\
&\begin{cases}
x_1^t &= J_{\frac{2\gamma}{w(n-1)}A_1}\left(\frac{2\mu}{w(n-1)} a_1 z_1^t\right)\\
x_i^t &= J_{\frac{2\gamma}{w(n-1)}A_i}\big(\frac{2}{w(n-1)}\big(w\sum_{j=1}^{i-1} x_j^t + \mu a_i z_i^t - \mu\sum_{j=1}^{i-1} b_j z_j^t - \gamma y_i\big)\big),\quad i\in\{2,\dots, n-1\}\\
x_n^t &= J_{\frac{2\gamma}{w(n-1)}A_n}\left(\frac{2}{w(n-1)}\left(w\sum_{j=1}^{n-1} x_j^t - \mu\sum_{j=1}^{n-1} b_j z_j^t - \gamma y_n \right)\right)\\
z_i^{t+1} &= z_i^t - \lambda_t \mu(a_i x_i^t - b_i \sum_{j=i+1}^n x_j^t), \quad i\in\{1,\dots,n-1\},
\end{cases}
\end{align}
where for $i,j\in\{1,\dots,n-1\}$, $a_i = \sqrt{\frac{(n-i)n}{n-i+1}}$ and $b_j = \sqrt{\frac{n}{(n-j)(n-j+1)}}$. This algorithm becomes \cite[Equation~(49)]{DPTT26} without dual terms when $\bB=0$, but differs in the choices of $H$ and $K$.
\end{enumerate}
As discussed in Remark~\ref{r:compare_stepsize}, Assumption~\ref{a:PSD} is satisfied whenever $\gamma < \frac{1}{\lambda_{\max}(\Theta)}$ and $\alpha = \gamma \lambda_{\max}(\Theta)$, where $\Theta$ is defined in \eqref{def:Theta}.
\end{remark}

\begin{remark}[Forward-reflected-backward algorithm and its variant]
    We show that the sequential algorithm presented in Remark~\ref{r:algorithms} can recover the forward-reflected-backward and its variant as special cases. Let $n = 3$, $r=1$, $p=2$, and $w = \mu = 1$. The algorithm \eqref{algo:seq} then becomes 
    \begin{align*}
    \begin{cases}
        x_1^t &= J_{2\gamma A_1}(2z_1^t)\\
        x_2^t &= J_{\gamma A_2}(x_1^t + z_2^t - z_1^t - \gamma B_1 x_1^t - \gamma C_1 x_1^t)\\
        x_3^t &= J_{2\gamma A_3}(2(x_2^t - z_2^t + \gamma B_1 x_1^t - \gamma B_1 x_2^t - \gamma C_2 x_2^t))\\
        z_1^{t+1} &= z_1^t - \lambda_t (x_1^t - x_2^t)\\
        z_2^{t+1} &= z_2^t - \lambda_t (x_2^t - x_3^t).
    \end{cases}
    \end{align*}
    Let $A_1 = 0$, $A := A_2$, $A_3 = 0$, $B := B_1$, $C := C_1$, $C_2 = 0$, and $\lambda_t = \frac{1}{2}$. The problem becomes 
    \begin{align*}
        \text{find~} x \ \text{~such that~}\ 0 \in A x + B x + C x,
    \end{align*}
    and the algorithm reduces to 
    \begin{align*}
        \begin{cases}
            x_1^t &= 2z_1^t\\
            x_2^t &= J_{\gamma A}(z_1^t + z_2^t - \gamma B x_1^t - \gamma C x_1^t)\\
            x_3^t &= 2(x_2^t - z_2^t + \gamma B x_1^t - \gamma B x_2^t)\\
            z_1^{t+1} &= z_1^t - \frac{1}{2}(2z_1^t - x_2^t) = \frac{1}{2}x_2^t\\
            z_2^{t+1} &= z_2^t - \frac{1}{2}(x_2^t - x_3^t).
        \end{cases}
    \end{align*}
    This implies that $z_1^t = \frac{1}{2}x_2^{t-1}$, $x_1^t = x_2^{t-1}$, and
    \begin{align*}
        z_2^{t+1} &= z_2^t - \frac{1}{2}(x_2^t - 2(x_2^t - z_2^t + \gamma B x_1^t - \gamma B x_2^t))= \frac{1}{2}x_2^t + \gamma B x_2^{t-1} - \gamma B x_2^t.
    \end{align*}
    Altogether, we have
    \begin{align}
    \label{algo:variant_FRB}
        x_2^t = J_{\gamma A}(x_2^{t-1} - 2\gamma B x_2^{t-1} + \gamma B x_2^{t-2} - \gamma C x_2^{t-1}),
    \end{align}
    which is precisely the variant of the forward-reflected-backward algorithm studied in \cite[Theorem~5.2]{MT20}. By Remark~\ref{r:compute_stepsize} and the given coefficient matrices, we have $\Omega = 0$ and $\gamma \leq \frac{\alpha}{\ell + \frac{1}{2}\beta}$. Note that to choose $\lambda_t = \frac{1}{2}$, we require $0 < \alpha < \frac{1}{2}$, thus $\gamma < \frac{1}{2\ell + \beta}$. This condition is slightly more restrictive than that of the variant of the forward-reflected-backward algorithm, namely $\gamma <\frac{1}{2\ell + \frac{1}{2}\beta}$. However, it is less restrictive than the bound obtained by applying \cite[Theorem 2.5]{MT20} directly, which requires $\gamma <\frac{1}{2\ell + 2\beta}$. Finally, when $C=0$, the algorithm \eqref{algo:variant_FRB} reduces to the forward-reflected-backward algorithm and recover the same stepsize bound $\gamma < \frac{1}{2\ell}$, as that established in \cite[Theorem 2.5]{MT20}. This also highlights the benefit of introducing the parameter $\alpha$ in Assumption~\ref{a:PSD}, which allows for a less restrictive bound on $\gamma$ than the condition in \emph{Case 2} of \cite[Remark~4.2]{DTT26}, where only $\gamma < \frac{1}{4\ell}$ is guaranteed.
\end{remark}

For each $i\in \{1, \dots, n\}$, let $U_i$ be a closed linear subspace of $\mathcal{H}$ and let $A_i =N_{U_i}$, where
\begin{align*}
N_{U_i}x=
\begin{cases}
U_i^\perp &\text{if~} x\in U_i, \\
\varnothing &\text{otherwise}.
\end{cases}    
\end{align*}
We set 
\begin{align*}
U :=U_1\cap \dots \cap U_n
\text{~~and~~}
\mathcal{U} :=U_1\times \dots \times U_n \subseteq \mathcal{H}^n.  
\end{align*}

\begin{theorem}[Fixed-point decomposition and explicit limits]
\label{t:explicit}
Suppose that, for each $i\in \{1, \dots, n\}$, $A_i =N_{U_i}$, where $U_i$ is a closed linear subspace, for each $k\in \{1, \dots, r\}$ and each $j\in \{1, \dots, p\}$, $B_k$ and $C_j$ are linear operators, that $Z\neq \varnothing$, and that Assumptions~\ref{a:stand}, \ref{a:explicit}, and \ref{a:PSD} hold. Let $(\bx^t)_{t\in\mathbb{N}}$ and $(\bz^t)_{t\in\mathbb{N}}$ be the sequences generated by Algorithm~\ref{algo:full} with an initial point $\bz^0\in \mathcal{H}^{m}$ and $(\lambda_t)_{t\in\mathbb{N}}$ in $[0,1-\alpha]$ satisfying $\sum_{t=1}^{+\infty} \lambda_t(1-\frac{\lambda_t}{1-\alpha}) = +\infty$. Then the following hold:
\begin{enumerate}
\item\label{t:explicit_ker}
$Z$ is a linear subspace and $Z\subseteq W := U\cap \bigcap_{j=1}^p \ker C_j\cap \bigcap_{k=1}^r \{x\in \mathcal{H}:\langle x, B_kx\rangle =0\}$, where equality holds if
\begin{align}\label{eq:compatibility-W}
\forall x\in W,\quad P\bB(\bOne x)\in \mathcal{U}^\perp.
\end{align}

\item\label{t:explicit_decomp}
$E :=\ker S$ is a closed linear subspace and, for each $x\in Z$, there exists a unique $\Lambda x\in E^\perp$ such that $S(\Lambda x) =\bOne x$. The mapping $\Lambda\colon Z\to E^\perp$ is linear and injective, $\Lambda Z :=\{\Lambda x: x\in Z\}$ is a closed linear subspace, and
\begin{align}\label{eq:FixT_decomp}
\Fix T =\Lambda Z\oplus E.
\end{align}
Consequently, for all $\bz\in \mathcal{H}^m$, $P_{\Fix T}\bz =P_{\Lambda Z}\bz +P_E\bz$.

\item\label{t:explicit_limit}
$\bz^t\to \bar{\bz} =P_{\Lambda Z}\bz^0 +P_E\bz^0\in \Fix T$ and $\bx^t\to \bOne\bar{x}$, where $\bar{x}$ is the unique point of $Z$ satisfying $\Lambda\bar{x} =P_{\Lambda Z}\bz^0$.

\item\label{t:explicit_compatible}
Suppose that
\begin{align}\label{eq:compatibility}
\forall x\in Z,\quad P\bB(\bOne x)\in \mathcal{U}^\perp.
\end{align}
Then $Z$ is a closed linear subspace and, for all $x\in Z$, $\Lambda x =\boldsymbol{\eta}x$, 
Consequently, $\Lambda Z=\boldsymbol{\eta}Z$ and, for all $\bz\in \mathcal{H}^m$, $P_{\Lambda Z}\bz =\frac{\boldsymbol{\eta}}{\|\boldsymbol{\eta}\|^2} P_Z(\boldsymbol{\eta}^\top \bz)$. In particular, the limits in \ref{t:explicit_limit} are explicitly given by $\bar{\bz} =\frac{\boldsymbol{\eta}}{\|\boldsymbol{\eta}\|^2} P_Z(\boldsymbol{\eta}^\top \bz^0) +P_E\bz^0$ and $\bar{x} =\frac{1}{\|\boldsymbol{\eta}\|^2} P_Z(\boldsymbol{\eta}^\top \bz^0)$.
\end{enumerate}
\end{theorem}
\begin{proof}
\ref{t:explicit_ker} \& \ref{t:explicit_decomp}: Let $x\in Z$, i.e., $0\in \sum_{i=1}^n N_{U_i}(x) +\sum_{k=1}^r B_kx +\sum_{j=1}^p C_jx$. We must have that, for all $i\in \{1, \dots, n\}$, $x\in U_i$, so $x\in U$ and 
\begin{align*}
\sum_{k=1}^r B_kx +\sum_{j=1}^p C_jx\in -\sum_{i=1}^n U_i^\perp =\sum_{i=1}^n U_i^\perp \subseteq U^\perp.
\end{align*}
It follows that
\begin{align*}
\sum_{k=1}^r \langle x, B_kx\rangle +\sum_{j=1}^p\langle x, C_jx\rangle =0.
\end{align*}
For each $k\in \{1, \dots, r\}$ and each $j\in \{1, \dots, p\}$, since $B_k$ is monotone and linear, and $C_j$ is cocoercive and linear,
\begin{align*}
\langle x, B_kx\rangle \geq 0
\text{~~and~~}
\langle x, C_jx\rangle\geq \frac{1}{\beta_j}\|C_jx\|^2 \geq 0.
\end{align*}
Therefore, $\langle x, B_k x\rangle =0$ and $C_jx =0$, which implies $Z \subseteq W$. 

Assume that \eqref{eq:compatibility-W} holds. Let $x\in W$ and write $P\bB(\bOne x) =(b_1, \dots, b_n)$. We have, for all $i\in \{1, \dots, n\}$ and all $j\in \{1, \dots, p\}$ that  $x\in U_i$, $C_jx =0$, and, by~\eqref{eq:compatibility-W}, $-b_i\in U_i^\perp =N_{U_i}(x)$. Noting that $\sum_{i=1}^n b_i =\bOne^\top P\bB(\bOne x) =\sum_{k=1}^r B_kx$, we obtain
\begin{align*}
0\in \sum_{i=1}^n N_{U_i}(x) +\sum_{k=1}^r B_kx +\sum_{j=1}^p C_jx,
\end{align*}
which means $x\in Z$. Therefore, $W\subseteq Z\subseteq W$, and so $Z = W$.

Since $S$ is linear and continuous by Theorem~\ref{t:cvg_strong}\ref{t:ST}, both $E =\ker S$ and $\Fix T =\ker(\Id -T) =\ker (M^\top S)$ are closed linear subspaces. Since $E =\ker S\subseteq \ker(M^\top S) =\Fix T$, for each $\bz\in \Fix T$, we have $P_E\bz\in E\subseteq \Fix T$, and so $\bz -P_E\bz\in \Fix T\cap E^\perp$. Therefore,
\begin{align}\label{eq:FixT-E-decomp}
\Fix T =(\Fix T\cap E^\perp)\oplus E.    
\end{align}
By Lemma~\ref{l:fixed_points}, $S(\Fix T) =\{\bOne x: x\in Z\}$ is a linear subspace, and the linearity and injectivity of $x\mapsto \bOne x$ imply that $Z$ is a linear subspace.

Let $x\in Z$. By Lemma~\ref{l:fixed_points}, there exists $\bz\in \Fix T$ such that $S\bz =\bOne x$. Since $P_E\bz\in E =\ker S$,
\begin{align*}
\by :=\bz -P_E\bz\in E^\perp \text{~~and~~} S\by =S\bz =\bOne x.
\end{align*}
If there is also $\bz'\in E^\perp$ satisfying $S\by' =\bOne x$, then $\by' -\by\in \ker S\cap E^\perp =\{0\}$. Hence, $\Lambda$ is well defined.

Next, for $x, y\in Z$ and $\alpha, \beta\in \mathbb{R}$, we have $\alpha\Lambda x +\beta\Lambda y\in E^\perp$ and, since $S$ is linear, $S(\alpha\Lambda x +\beta\Lambda y) =\bOne(\alpha x +\beta y)$. By uniqueness, $\Lambda(\alpha x +\beta y) =\alpha\Lambda x +\beta\Lambda y$. Thus, $\Lambda$ is linear. Moreover, $\Lambda x=0$ implies $\bOne x=S(\Lambda x)=0$, so $\Lambda$ is injective.

Now, it follows directly from the definition of $\Lambda$ that $\Lambda Z\subseteq \Fix T\cap E^\perp$. Conversely, if $\bz\in \Fix T\cap E^\perp$, then, by Lemma~\ref{l:fixed_points}, there exists $x\in Z$ such that $S\bz =\bOne x$, which yields $\bz= \Lambda x\in \Lambda Z$. We deduce that $\Lambda Z =\Fix T\cap E^\perp$, which together with \eqref{eq:FixT-E-decomp} implies \eqref{eq:FixT_decomp}. 

\ref{t:explicit_limit}: Since $A_i =N_{U_i}$ is a linear relation, Theorem~\ref{t:cvg_strong} implies that $\bz^t\to \bar{\bz} =P_{\Fix T}(\bz^0)$. In view of \ref{t:explicit_decomp}, $\bar{\bz} =P_{\Lambda Z}\bz^0 +P_E\bz^0$. Since $P_{\Lambda Z}\bz^0\in \Lambda Z$ and $\Lambda$ is injective, there is a unique $\bar{x}\in Z$ satisfying $\Lambda\bar{x} =P_{\Lambda Z}\bz^0$. Since $P_E\bz^0\in \ker S$, the continuity of $S$ gives
\begin{align*}
\bx^t =S\bz^t\to S(P_{\Lambda Z}\bz^0 +P_E\bz^0) =S(\Lambda\bar{x}) =\bOne\bar{x}.
\end{align*}

\ref{t:explicit_compatible}: For all $x\in Z$, we have from \ref{t:explicit_ker} that $x\in U\cap \bigcap_{j=1}^p \ker C_j$, which yields $\bA(\bOne x) = \mathcal{U}^\perp$ and $\bC(\bOne x) =0$, so
\begin{align}
\label{eq:U_perp}
(\bA +P\bB +H\bC)(\bOne x) =\mathcal U^\perp +P\bB(\bOne x) =\mathcal{U}^\perp,
\end{align}
where the last equality follows from \eqref{eq:compatibility}.  

We next show that
\begin{align}\label{eq:muU-E}
\boldsymbol{\eta}U\subseteq E^\perp.
\end{align}
Since $\boldsymbol{\eta} =M^\dagger(D -N)\bOne\in \ran M^\top$, there exists $\bv\in \mathbb{R}^n$ such that $\boldsymbol{\eta} =M^\top\bv$. Let $\be\in E \subseteq \Fix T$. Then $S\be =0 =\bOne 0$, which together with Lemma~\ref{l:fixed_points} and \eqref{eq:U_perp} implies that $M\be\in \mathcal{U}^\perp$, i.e., for all $i\in \{1, \dots, n\}$, $(M\be)_i\in U_i^\perp$. Therefore, for all $x \in U \subseteq U_i$,
\begin{align*}
\langle \boldsymbol{\eta}x, \be\rangle =\langle M^\top\bv x, \be\rangle =\langle\bv x, M\be\rangle =\sum_{i=1}^n v_i\langle x, (M\be)_i\rangle =0,
\end{align*}
which proves \eqref{eq:muU-E}.

Now, let $x\in Z$. By \eqref{eq:U_perp}, $M(\boldsymbol{\eta}x -\boldsymbol{\eta}x) = 0 \in \mathcal{U}^\perp$, which implies $S(\boldsymbol{\eta}x) =\bOne x$ due to Lemma~\ref{l:fixed_points}. Moreover, $Z\subseteq U$ by \ref{t:explicit_ker}, so \eqref{eq:muU-E} gives $\boldsymbol{\eta}x\in E^\perp$. By the uniqueness of $\Lambda x$, we obtain $\Lambda x =\boldsymbol{\eta}x$. Thus, $\Lambda Z =\boldsymbol{\eta}Z$.

Since $\Lambda Z$ is closed and $\boldsymbol{\eta}\neq 0$, $Z$ is closed. Let $\bz\in \mathcal H^m$ be arbitrary and set $x_{\bz} :=\frac{1}{\|\boldsymbol{\eta}\|^2}P_Z(\boldsymbol{\eta}^\top \bz)\in Z$. For all $y\in Z$,
\begin{align*}
\langle \bz -\boldsymbol{\eta}x_{\bz}, \boldsymbol{\eta}y\rangle =\langle \boldsymbol{\eta}^\top \bz, y\rangle -\|\boldsymbol{\eta}\|^2\langle x_{\bz}, y\rangle =\left\langle\boldsymbol{\eta}^\top\bz-P_Z(\boldsymbol{\eta}^\top\bz),y\right\rangle =0.
\end{align*}
As $\boldsymbol{\eta}x_{\bz}\in \boldsymbol{\eta}Z =\Lambda Z$, it follows that
\begin{align*}
P_{\Lambda Z}\bz =P_{\boldsymbol{\eta}Z}\bz =\frac{\boldsymbol{\eta}}{\|\boldsymbol{\eta}\|^2} P_Z(\boldsymbol{\eta}^\top \bz).
\end{align*}
In particular, $P_{\Lambda Z}\bz^0 =\frac{\boldsymbol{\eta}}{\|\boldsymbol{\eta}\|^2} P_Z(\boldsymbol{\eta}^\top \bz^0)$ and $\bar{\bz} =\frac{\boldsymbol{\eta}}{\|\boldsymbol{\eta}\|^2} P_Z(\boldsymbol{\eta}^\top \bz^0) +P_E\bz^0$. By \ref{t:explicit_limit},
\begin{align*}
\boldsymbol{\eta}\bar{x} =\Lambda\bar{x} =P_{\Lambda Z}\bz^0 =\frac{\boldsymbol{\eta}}{\|\boldsymbol{\eta}\|^2} P_Z(\boldsymbol{\eta}^\top \bz^0).
\end{align*}
Since $\boldsymbol{\eta}\neq 0$, the mapping $x\mapsto \boldsymbol{\eta}x$ is injective, and so $\bar{x} =\frac{1}{\|\boldsymbol{\eta}\|^2} P_Z(\boldsymbol{\eta}^\top \bz^0)$.
\end{proof}

\begin{remark}[Limit points of some algorithms]
We illustrate how to derive the limit point for each instance of Algorithm~\ref{algo:full}. In view of Theorem~\ref{t:explicit}\ref{t:explicit_compatible}, we have $\bar{x} =P_Z(\boldsymbol{\eta}^\dagger\bz^0)$. Therefore, given an initial point $\bz^0$, the limit point $\bar{x}$ is obtained as soon as $\boldsymbol{\eta}= (\eta_1,\dots,\eta_{n-1}) =M^\dagger (D -N)\bOne$ is specified. The values of $\boldsymbol{\eta}$ and limit points of the sequences generated by six algorithms are provided in Table~\ref{tab:my-table}. For the case when $Z=U$, details about calculations of $\boldsymbol{\eta}$ and $\bar{x}$ can be found in \cite[Section~4.4 and 4.6]{ABCL26}.
\begin{table}[!htbp]
\renewcommand{\arraystretch}{1.5}
\centering
\begin{tabular}{|c|c|c|}
\cline{1-3}
Algorithm    & $\eta_j$          & Limit point $\bar{x}$  \\ \cline{1-3}
Complete-parallel down & $\frac{n+1 -2j}{2}w\mu$ & $\frac{12}{w^2\mu^2(n-1)((n-1)^2 + 2)} P_Z (\sum_{j=1}^{n-1} \eta_j z_j)$ \\ \cline{1-3}
Ring-sequential  & $w\mu$        & $\frac{1}{w\mu(n-1)}P_Z (\sum_{j=1}^{n-1} z_j)$ \\ \cline{1-3}
Parallel up     & $\frac{1}{2}w\mu$       & $\frac{2}{w\mu(n-1)}P_Z (\sum_{j=1}^{n-1} z_j)$ \\ \cline{1-3}
Parallel down     & $\frac{1}{2}w\mu$     & $\frac{2}{w\mu(n-1)}P_Z (\sum_{j=1}^{n-1} z_j)$ \\ \cline{1-3}
Sequential       & $\frac{1}{2}w\mu$      & $\frac{2}{w\mu(n-1)}P_Z (\sum_{j=1}^{n-1} z_j)$ \\ \cline{1-3}
Complete     & $w\mu\sqrt{\frac{(n-j)(n-j+1)}{4n}}$          & $\frac{12}{w^2\mu^2(n^2-1)}P_Z (\sum_{j=1}^{n-1} \eta_j z_j)$ \\ \cline{1-3}
\end{tabular}
\caption{Limit points of the sequences generated by algorithms in Remark~\ref{r:algorithms}.}
\label{tab:my-table}
\end{table}
\end{remark}

\begin{remark}[An explicit form of $Z$]
\label{r:compatibility}
Regarding the assumptions in Theorem~\ref{t:explicit}\ref{t:explicit_ker}\&\ref{t:explicit_compatible}, we note that \eqref{eq:compatibility-W} implies \eqref{eq:compatibility}. A sufficient condition for \eqref{eq:compatibility-W} is
\begin{align}\label{eq:kercondition}
U\cap \left(\bigcap_{k=1}^r \{x\in \mathcal{H}:\langle x, B_kx\rangle =0\}\right) \subseteq \bigcap_{k=1}^r \ker B_k,
\end{align}
in which case, 
\begin{align*}
Z =W =U\cap \bigcap_{k=1}^r \ker B_k\cap \bigcap_{j=1}^p \ker C_j.    
\end{align*}
Indeed, it is clear that $U\cap \bigcap_{k=1}^r \ker B_k\cap \bigcap_{j=1}^p \ker C_j\subseteq W$, while the reverse inclusion follows from \eqref{eq:kercondition}. Hence, $W =U\cap \bigcap_{k=1}^r \ker B_k\cap \bigcap_{j=1}^p \ker C_j$. In particular, for all $x\in W$, we have $x\in \bigcap_{k=1}^r \ker B_k$, so \eqref{eq:compatibility-W} holds, and Theorem~\ref{t:explicit}\ref{t:explicit_ker} yields $Z =W$.

It is worth noting that \eqref{eq:kercondition} is satisfied if there exists $\ell >0$ such that
\begin{align*}
\forall x\in U, \quad \sum_{k=1}^r \langle x, B_kx\rangle\geq \ell\sum_{k=1}^r \|B_kx\|^2.    
\end{align*}
\end{remark}

The following results are the fixed-point decomposition and explicit limit points for the case where $\bB=0$. Furthermore, when $\bC=0$, we recover the results in \cite[Theorem~4.8]{ABCL26} which require $M$ to be injective and $m=n-1$.
\begin{corollary}
Suppose that, for each $i\in \{1, \dots, n\}$, $A_i =N_{U_i}$, where $U_i$ is a closed linear subspace, for each $j\in\{1,\dots,p\}$, $C_j$ is a linear operator, that $Z = \zer\left(\sum_{i=1}^n N_{U_i} + \sum_{j=1}^p C_j\right)\neq \varnothing$, and that Assumptions~\ref{a:stand}, \ref{a:explicit}, and \ref{a:PSD} hold. Let $(\bx^t)_{t\in\mathbb{N}}$ and $(\bz^t)_{t\in\mathbb{N}}$ be the sequences generated by Algorithm~\ref{algo:full} with an initial point $\bz^0\in \mathcal{H}^{m}$ and $(\lambda_t)_{t\in\mathbb{N}}$ in $[0,1-\alpha]$ satisfying $\sum_{t=1}^{+\infty} \lambda_t(1-\frac{\lambda_t}{1-\alpha}) = +\infty$. Then the following hold:
\begin{enumerate}
\item
$Z =U\cap \bigcap_{j=1}^p \ker C_j$, $E =\ker S$ are closed linear subspaces, and $\Fix T = \boldsymbol{\eta} Z\oplus E$.
\item
$\bz^t\to \bar{\bz} =P_{\Fix T}(\bz^0) =\frac{\boldsymbol{\eta}}{\|\boldsymbol{\eta}\|^2} P_Z(\boldsymbol{\eta}^\top \bz^0) +P_E\bz^0$ and $\bx^t\to \bOne\bar{x}$, where $\bar{x} =\frac{1}{\|\boldsymbol{\eta}\|^2} P_Z(\boldsymbol{\eta}^\top \bz^0)$.
\end{enumerate}    
\end{corollary}
\begin{proof}
In this case, $\bB = 0$, so \eqref{eq:compatibility-W} and \eqref{eq:compatibility} are automatically satisfied. It follows from Theorem~\ref{t:explicit}\ref{t:explicit_ker}\&\ref{t:explicit_compatible} that $Z =U\cap \bigcap_{j=1}^p \ker C_j$ is a closed linear subspace and $\Lambda Z = \boldsymbol{\eta}Z$. The conclusions then follow from Theorem~\ref{t:explicit}\ref{t:explicit_decomp}--\ref{t:explicit_compatible}.
\end{proof}

\section{Numerical experiments}
\label{s:numerics}
In this section, we provide a numerical experiment to demonstrate the application of Algorithm~\ref{algo:full} to a particular instance of problem \eqref{prob:ABC}. We consider the quadratic minimization problem with multiple-set split feasibility constraints
\begin{align}
\label{prob:gen}
    \min_{x\in\mathbb{R}^d} \quad \sum_{j=1}^p \left(\frac{1}{2} x^\top S_j x\right) \quad\text{~s.t.~} x \in \mathcal{C} := \bigcap_{i=1}^n \mathcal{C}_i \ \text{~and~} \ L x \in \mathcal{Q} := \bigcap_{k=1}^r \mathcal{Q}_k,
\end{align}
where for $i\in\{1,\dots,n\}$, $k\in\{1,\dots,r\}$, $\mathcal{C}_i \subseteq \mathbb{R}^d$ and $\mathcal{Q}_k \subseteq \mathbb{R}^s$ are nonempty closed convex sets, $L \in \mathbb{R}^{s \times d}$, and, for each $j\in\{1,\dots,p\}$, $S_j \in \mathbb{R}^{d \times d}$ is symmetric and positive semidefinite. When for all $j$, $S_j=0$, this problem reduces to the \emph{multiple-sets split feasibility problem} \cite{CEKB05}. Using the Fenchel representation of the indicator function of each $\mathcal{Q}_k$, problem \eqref{prob:gen} can be expressed as
\begin{align}
    \min_{x\in\mathbb{R}^d} \max_{y_1,\dots,y_r \in \mathbb{R}^s} \quad \sum_{j=1}^p \left(\frac{1}{2} x^\top S_j x\right) + \iota_{\cap_{i=1}^n \mathcal{C}_i}(x) + \sum_{k=1}^r \scal{y_k}{L x} - \sum_{k=1}^r \iota_{\mathcal{Q}_k}^* (y_k),
\end{align}
where $\iota_{\mathcal{Q}_k}^*= \sigma_{\mathcal{Q}_k}$ is the support function of $\mathcal{Q}_k$. To derive optimality conditions involving the normal cones of the individual constraint sets, we impose the following qualification conditions \cite[Corollary~16.50(v)]{BC17}: $\bigcap_{i=1}^n \operatorname{ri} \mathcal{C}_i \neq \varnothing$ and $\bigcap_{k=1}^r \operatorname{ri} \mathcal{Q}_k \neq \varnothing$, where $\operatorname{ri}$ denotes the relative interior. Then the optimality conditions are
\begin{align}
\label{num:opt_con}
    0 \in \sum_{j=1}^p S_j x + \sum_{i=1}^n N_{\mathcal{C}_i}(x) + \sum_{k=1}^r L^* y_k, \quad 0 \in N_{\mathcal{Q}_k}^{-1}(y_k) - L x, \quad k\in\{1,\dots,r\},
\end{align}
where $N_{\mathcal{Q}_k}^{-1} = \partial \iota_{\mathcal{Q}_k}^*$, followed from \cite[Theorem~16.29(i)\&(iv)]{BC17} and the fact the $N_{\mathcal{Q}_k} = \partial\iota_{\mathcal{Q}_k}$. Let $A_1 = \diag(N_{\mathcal{C}_1}, N_{\mathcal{Q}_1}^{-1}, \dots, N_{\mathcal{Q}_r}^{-1})$, $A_i = \diag(N_{\mathcal{C}_i}, 0, \dots, 0)$ for $i\in\{2,\dots,n\}$. Let
$$(B_k)_{ij} = 
\begin{cases}
    L^* \quad &\text{if~} i=1, j=k+1\\
    -L \quad &\text{if~} i=k+1, j=1,\\
    0 \quad &\text{otherwise},
\end{cases}$$
and let $C_j = \diag(S_j, 0, \dots, 0)$ for $j\in\{1,\dots,p\}$. This problem is equivalent to finding $(x,y_1,\dots,y_r) \in \mathbb{R}^d \times \mathbb{R}^{s} \times \dots \times  \mathbb{R}^s$ such that
\begin{align}
    0 \in \sum_{i=1}^n A_i(x,y_1,\dots,y_r)  + \sum_{k=1}^r B_k (x,y_1,\dots,y_r) + \sum_{j=1}^p C_j(x,y_1,\dots,y_r).
\end{align}
To show how limit points of the sequences generated by several popular algorithms can be computed explicitly, we consider the case where $\mathcal{C}_i$ and $\mathcal{Q}_k$ are closed linear subspaces. In this case, $N_{\mathcal{Q}_k}^{-1} = N_{\mathcal{Q}_k^\perp}$, thus $A_1=\diag(N_{\mathcal{C}_1}, N_{\mathcal{Q}_1^\perp}, \dots, N_{\mathcal{Q}_r^\perp})$. Moreover, $A_i$ is maximally monotone, $B_k$ is monotone and $\|L\|_2$-Lipschitz continuous, and $C_j$ is $\frac{1}{\|S_j\|_2}$-cocoercive. In addition, $A_i$ is a linear relation, as it is defined by the normal cones to closed linear subspaces and theirs orthogonal complements, whereas $B_k$ and $C_j$ are linear operators.

\paragraph{Data generation.} Let $n \geq 3$, $r, p \geq 1$, $d\geq 4$ and $s \geq 2$. Let $\{v_1,\dots,v_d\}$ be an orthonormal basis of $\mathbb{R}^d$. The data generation process is designed so that the constraints define a 2-dimensional plane $\lspan\{v_1, v_2\}$, while the optimal solution set for the problem is $\lspan\{v_1\}$. This construction guarantees that the qualification conditions hold and the problem admits nonunique and nontrivial solutions. Specifically, the procedure works as follows:
\begin{itemize}
    \item
    Defining $S_j$: For each $j\in\{1,\dots,p\}$, define $S_j = \frac{j}{j+1}(\Id - v_1 v_1^\top - v_4 v_4^\top)$. Each $S_j$ is symmetric and positive semidefinite with $\ker S_j = \lspan\{v_1, v_4\}$. Consequently, $\bigcap_{j=1}^p \ker S_j=\lspan\{v_1, v_4\}$.
    \item
    Defining $\mathcal{C}_i$: We construct the constraint subspace $\mathcal{C}_i$ so that their intersection is exactly a space spanned by $\{v_1, v_2, v_3\}$.  and define $\mathcal{C}_1 = \lspan\{v_1, v_2,v_3\}$. For each $i\in\{2,\dots,n\}$, choose a distinct nonempty subset $I_i \subseteq \{v_4,\dots,v_d\}$ and define $\mathcal{C}_i = \lspan(\{v_1,v_2,v_3\} \cup I_i)$. Each $\mathcal{C}_i$ is a closed linear subspace and their intersection is $\mathcal{C} = \lspan\{v_1,v_2,v_3\}$. Since there are $2^{d-3}-1$ such subsets, the number of subspaces must satisfy $n \leq 2^{d-3}$.
    \item
    Defining $\mathcal{Q}_k$: We construct closed linear subspaces whose intersection is a line spanned by $q_1$. First, choose an orthonormal basis $\{q_1,\dots,q_s\}$ of $\mathbb{R}^s$ and define $\mathcal{Q}_1 = \lspan\{q_1\}$. For each $k\in\{2,\dots,r\}$, select a distinct nonempty subset $I_k \subseteq \{q_2,\dots,q_s\}$ that must contain $q_2$ and define $\mathcal{Q}_k =\lspan(\{q_1\} \cup I_k\})$. Each $\mathcal{Q}_k$ is a closed linear subspace and their intersection is $\mathcal{Q} = \lspan\{q_1\}$. Since there are $2^{s-2}$ such subsets, the number of subspaces must satisfy $r\leq 2^{s-2} + 1$.
    \item
    Defining $L$: Let $L = a q_1 v_2^\top + b q_2 v_3^\top \in \mathbb{R}^{s \times d}$, where $a\in \mathbb{R}$ and $b\in\mathbb{R}\setminus \{0\}$ are random numbers. For any $x \in \mathcal{C}$, we have $x = \alpha v_1 + \beta v_2 + \zeta v_3$ with $\alpha, \beta, \zeta \in \mathbb{R}$. By orthonormality, $Lx= a \beta q_1 + b \zeta q_2$. Since $q_1$ and $q_2$ are orthonormal, $\mathcal{Q}=\lspan\{q_1\}$, and $b\neq 0$, we have $Lx \in \mathcal{Q}$ if and only if $\zeta = 0$. Consequently, $\mathcal{C} \cap L^{-1}(\mathcal{Q}) = \lspan\{v_1,v_2\}$ where $L^{-1}(\mathcal{Q}) := \{ x \in \mathbb{R}^d : Lx \in \mathcal{Q}\}$.
\end{itemize}
By construction, the solution set of problem \eqref{prob:gen} is $\mathcal{C} \cap L^{-1}(\mathcal{Q}) \cap \bigcap_{j=1}^p \ker S_j = \lspan\{v_1\}$. Since $N_{\mathcal{C}_i}(x) = \mathcal{C}_i^\perp$ as $x \in \mathcal{C}\subseteq \mathcal{C}_i$ and $N_{\mathcal{Q}_k}(Lx) = \mathcal{Q}_k^\perp$ as $Lx \in \mathcal{Q}\subseteq \mathcal{Q}_k$, the optimality conditions \eqref{num:opt_con} become
\begin{align}
\label{num:verify}
    0 \in \sum_{i=1}^n \mathcal{C}_i^\perp + \sum_{k=1}^r L^* y_k, \quad y_k \in \mathcal{Q}_k^\perp, \quad k \in \{1,\dots, r\}.
\end{align}
By construction, $\mathcal{C}=\lspan\{v_1,v_2,v_3\}$. Using $\mathcal{C} = \bigcap_{i=1}^n \mathcal{C}_i$, we obtain $\sum_{i=1}^n \mathcal{C}_i^\perp = \mathcal{C}^\perp = \lspan\{v_4,\dots,v_d\}$. Moreover, for $y_k \in \mathcal{Q}_k^\perp$, by the definition of $L^*$, we have 
\begin{align*}
    \sum_{k=1}^r L^* y_k = a v_2  \langle q_1, \sum_{k=1}^r y_k\rangle + b v_3 \langle q_2, \sum_{k=1}^r y_k\rangle.
\end{align*}
Since $y_1 \in \mathcal{Q}_1^\perp = \lspan\{q_2,\dots, q_s\}$ and $\sum_{k=2}^r y_k \in \sum_{k=2}^r \mathcal{Q}_k^\perp \subseteq \lspan\{q_3,\dots,q_s\}$, the coefficient of $v_2$ vanishes, while the inner product involving $q_2$ reduces to $\langle q_2, y_1\rangle$. Hence, $\sum_{k=1}^r L^* y_k = b v_3 \langle q_2, y_1 \rangle$. The first term in \eqref{num:verify} therefore becomes $0 \in \sum_{i=1}^n \mathcal{C}_i^\perp + \sum_{k=1}^r L^* y_k = \mathcal{C}^\perp + b v_3 \langle q_2, y_1 \rangle$. Since $v_3 \perp \mathcal{C}^\perp$ and $b\neq 0$, this condition holds if and only if $\langle q_2, y_1\rangle = 0$, which means $y_1 \in \mathcal{Q}_1^\perp \cap \lspan\{q_2\}^\perp = \lspan\{q_3,\dots,q_s\}$. Combining this with the primal constraints, the primal-dual optimal solution set is
\begin{align}
\label{num:Z}
    Z = \lspan\{v_1\} \times \lspan\{q_3,\dots,q_s\} \times \mathcal{Q}_2^\perp \times \dots \times \mathcal{Q}_r^\perp.
\end{align}

To apply Theorem~\ref{t:explicit}\ref{t:explicit_compatible}, we first verify the condition \eqref{eq:compatibility}. Specifically, for $i\in\{1,\dots,n\}$ and $(x,y_1,\dots,y_r)\in Z$, 
\begin{align*}
    \sum_{k=1}^r P_{ik} B_k(x,y_1,\dots,y_r) \in A_i (x,y_1,\dots,y_r)= U_i^\perp.
\end{align*}
Using \eqref{num:Z}, we have $x\in\lspan\{v_1\}$, which implies $Lx = 0$. Furthermore, $y_1\in \lspan\{q_3,\dots,q_s\}$ and for $k\in\{2,\dots,r\}$, $y_k \in \mathcal{Q}_k^\perp \subseteq \sum_{l=2}^r \mathcal{Q}_l^\perp \subseteq \lspan\{q_3,\dots,q_s\}$. It follows from the definition of $L^*$ that for $k\in\{1,\dots,r\}$, $L^* y_k=0$. Consequently, by the definition of $B_k$, for $k\in\{1,\dots,r\}$, $B_k(x,y_1,\dots,y_r) =0$. Thus, for any matrix $P$ and $(x,y_1,\dots,y_r) \in Z$, $P\bB(\bOne (x,y_1,\dots,y_r)) =0 \in \mathcal{U}^\perp$, which satisfies condition \eqref{eq:compatibility}.

\begin{figure}[htbp]
\captionsetup{skip=0pt}
    \centering
    \begin{subfigure}{0.48\textwidth}
        \includegraphics[width=\linewidth]{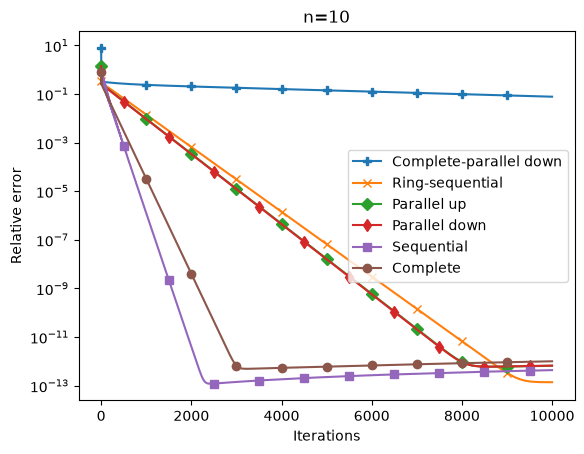}
    \end{subfigure}
    \begin{subfigure}{0.48\textwidth}
        \includegraphics[width=\linewidth]{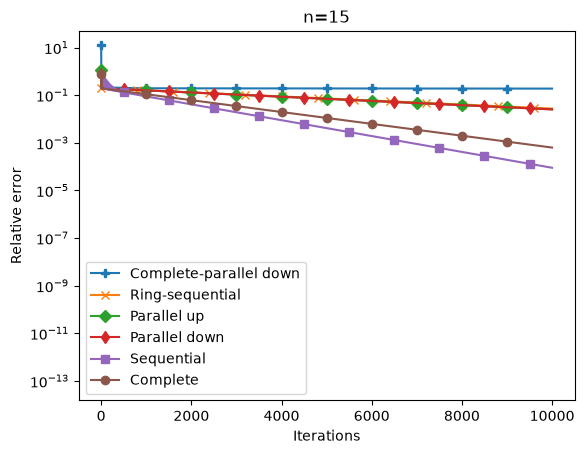}
    \end{subfigure}
    \caption{Comparing the performance of some algorithms}
    \label{fig:Lipschitz_comparison}
\end{figure}

We consider six algorithms in Remark~\ref{r:algorithms} with $w=\mu=1$ and their limit points are given in Table~\ref{tab:my-table}. To compare with results obtained by running the algorithms, we measure the performance by the relative error between the averaged iterate and the corresponding limit point, defined as 
\begin{align*}
    \text{relative error~} = \frac{\|(x, y_1,\dots,y_r)_{\rm avg}-(\bar{x}, \bar{y}_1,\dots,\bar{y}_r)\|}{\|(\bar{x}, \bar{y}_1,\dots,\bar{y}_r)\|},
\end{align*}
where $(x, y_1,\dots,y_r)_{\rm avg} = \frac{1}{n}(\sum_{i=1}^n x_i,\sum_{i=1}^n y_{1,i},\dots,\sum_{i=1}^n y_{r,i})$ and $(\bar{x}, \bar{y}_1,\dots,\bar{y}_r)$ denotes the corresponding limit point, computed as described in Table~\ref{tab:my-table}. The stepsize $\gamma$ and the relaxation parameter $\lambda_t$ are chosen empirically for all algorithms according to $\gamma = \hat{\gamma} \frac{1}{\lambda_{\max}(\Theta)}$ and $\lambda_t \equiv \lambda = \hat{\lambda}(1- \gamma \lambda_{\max}(\Theta))$, where $\Theta$ is given by \eqref{def:Theta}, and $(\hat{\gamma}, \hat{\lambda})= (0.5, 0.9)$ are scaled parameters. In the experiment, we set $d=s=10$. The results for the cases when $n = 10$ and $n = 15$ are given in Figure~\ref{fig:Lipschitz_comparison}.
The plots show consistent results in both cases, with the algorithms falling into three groups. The first group, comprising the sequential and complete algorithms, converge the fastest. Within this group, the sequential algorithm converges faster than the complete algorithm. The second group, consisting of the parallel up, parallel down, and ring-sequential algorithms, converges more slowly but still achieves high accuracy. The third group is the complete-parallel down algorithm which shows the lowest convergence with lower accuracy. The results also suggest that as $n$ increases, the performance of all algorithms decrease.

\paragraph{Acknowledgements.} The research of MND, MKT, and TDT was supported in part by Australian Research Council grant DP230101749.

\end{document}